\documentclass[11pt]{amsart}

\pdfoutput=1

\usepackage[text={420pt,650pt},centering]{geometry}

\usepackage{cancel}
\usepackage{esint,amssymb} 
\usepackage{graphicx}
\usepackage{MnSymbol}
\usepackage{mathtools} 
\usepackage[colorlinks=true, pdfstartview=FitV, linkcolor=blue, citecolor=blue, urlcolor=blue,pagebackref=false]{hyperref}
\usepackage{microtype}

\usepackage{bm}
\usepackage{dsfont}
\usepackage{mathrsfs}
\usepackage{xcolor}

\usepackage{wrapfig}
\usepackage{caption}

\newtheorem{proposition}{Proposition}
\newtheorem{theorem}[proposition]{Theorem}

\newtheorem{conjecture}[proposition]{Conjecture}
\theoremstyle{remark}

\theoremstyle{definition}

\numberwithin{equation}{section}
\numberwithin{proposition}{section}
\numberwithin{figure}{section}
\numberwithin{table}{section}

\newcommand{\R}{\mathbb{R}}

\newcommand{\E}{\mathbb{E}}

\newcommand{\EE}{\mathbf{E}}
\newcommand{\PP}{\mathbf{P}}

\newcommand{\ep}{\varepsilon}
\newcommand{\eps}{\varepsilon}

\renewcommand{\le}{\leqslant}
\renewcommand{\ge}{\geqslant}

\renewcommand{\subset}{\subseteq}

\newcommand{\td}{\widetilde}

\newcommand{\Ll}{\left}
\newcommand{\Rr}{\right}
\renewcommand{\d}{\mathrm{d}}
\newcommand{\dr}{\partial}

\newcommand{\1}{\mathbf{1}}

\newcommand{\mcl}{\mathcal}
\newcommand{\msf}{\mathsf}

\newcommand{\al}{\alpha}

\newcommand{\de}{\delta}
\newcommand{\si}{\sigma}

\newcommand{\la}{\left\langle}
\newcommand{\ra}{\right\rangle}

\newcommand{\bip}{{\mathrm{bip}}}

\newcommand{\mart}{\mathbf{Mart}}

\begin{document}

\author{Jean-Christophe Mourrat}
\address[Jean-Christophe Mourrat]{Department of Mathematics, ENS Lyon and CNRS, Lyon, France}

\keywords{}
\subjclass[2010]{}
\date{\today}

\title{Spin glasses and the Parisi formula}

\begin{abstract}
Spin glasses are models of statistical mechanics in which a large number of simple elements interact with one another in a disordered fashion. One of the fundamental results of the theory is the Parisi formula, which identifies the limit of the free energy of a large class of such models. Yet many interesting models remain out of reach of the classical theory, and direct generalizations of the Parisi formula yield invalid predictions. I will report here on some partial progress towards the resolution of this problem, which also brings a new perspective on classical results.
\end{abstract}

\maketitle

The aim of statistical mechanics is to describe the emergent properties of systems that are made of a large number of simple elements. Spin glasses are particular such models, in which there is a lot of ``disagreement'' between the elementary units of the system\footnote{Why such models may relate to glass is discussed in \cite{mourrat2024informal}.}.  Mathematically, this is usually modeled by introducing randomness into the interactions between the elements. 
 In the first section of this note, we make this concrete by presenting a basic spin glass called the Sherrington-Kirkpatrick (SK) model. We define the free energy of the model, and state a fundamental result, called the Parisi formula, that identifies the asymptotic behavior of the free energy in the limit of large system size. We also discuss surprising aspects of this formula, and in Section~\ref{s.uninverting}, we present a more recent alternative formulation. In Section~\ref{s.general}, some variants of the SK model are introduced for which the limit free energy is currently not known. This is the fundamental problem that has driven most of my work in the topic. In Section~\ref{s.pdes}, a point of view based on partial differential equations is introduced that allows us to formulate a natural conjecture for the limit free energy of these models. Partial results consistent with this conjecture are also presented. In Section~\ref{s.conjecture}, we discuss a promising connection between this point of view based on partial differential equations and the alternative representation of the Parisi formula that appeared in Section~\ref{s.uninverting}. The note ends with a short concluding section.

%
%
%
%
%
%
\section{The Parisi formula}
\label{s.parisi}

We start by introducing a basic spin glass called the Sherrington-Kirkpatrick (SK) model \cite{sherrington1975solvable}. We give ourselves independent Gaussian random variables $(W_{i,j})_{i,j \ge 1}$ of zero mean and unit variance, and for every $\sigma \in \R^N$, we set
\begin{equation}
\label{e.def.HN}
H_N(\sigma) := \frac 1 {\sqrt{N}} \sum_{i,j = 1}^N W_{i,j} \sigma_i \sigma_j.
\end{equation}
We are interested in questions such as: what is the limit of
\begin{equation}  
\label{e.def.max}
\frac 1 N \max_{\sigma \in \{-1,1\}^N} H_N(\sigma)
\end{equation}
as $N$ tends to infinity? This problem is often motivated with the following story. There are $N$ individuals $\{1,\ldots, N\}$ that need to be split into two groups. We can represent one such splitting using a vector $\sigma \in \{-1,1\}^N$, with the understanding that $\sigma_i$ indicates the group to which individual $i$ is assigned. The coefficients $W_{ij}$ represent how much individual $i$ likes individual $j$, and we would like to find an assignment $\sigma$ that maximizes global welfare, that is, a maximizer in $\{-1,1\}^N$ of
\begin{equation}  
\label{e.alt.max}
\sigma \mapsto \sum_{i,j = 1}^N W_{i,j} \1_{\{\sigma_i = \sigma_j\}}.
\end{equation}
For $\sigma \in \{-1,1\}^N$, we have $\sigma_i \sigma_j = 2\1_{\{\sigma_i = \sigma_j\}} - 1$, so the maximization of \eqref{e.alt.max} is essentially equivalent to that in \eqref{e.def.max}, up to an affine transformation. 

\begin{wrapfigure}{r}{0.3\linewidth}
    \centering
    \includegraphics[scale=1]{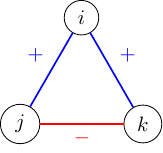}
    \captionsetup{width=0.95\linewidth}
\caption{{\small A simple situation with frustration. The coefficients $(W_{ij})$ suggest to set $\sigma_i = \sigma_j$, $\sigma_i = \sigma_k$, and $\sigma_j = - \sigma_k$, but we cannot realize these three conditions simultaneously.}}
    \label{f.frustration}
\end{wrapfigure}

Finding a configuration $\sigma \in \{-1,1\}^N$ that maximizes $H_N(\sigma)$ is a non-trivial task. Even with $N = 3$, we see in situations such as that depicted in Figure~\ref{f.frustration} that it will typically not be possible to find a configuration $\sigma$ such that for every $i$ and $j$, we have $(W_{i,j} + W_{j,i}) \sigma_i \sigma_j \ge 0$. In other words, certain pairs will be \emph{frustrated}: two individuals $i$ and $j$ may be assigned to different groups even though they would rather be together, or vice versa. The presence of these \emph{frustrations} is the key signature of spin glasses.

More fundamentally, one can show that the problem, given the coefficients $(W_{ij})$, of finding a configuration $\sigma \in \{\pm 1\}^N$ that maximizes $H_N$, is NP-hard in general. In fact, the problem is NP-hard even if we only aim to find a configuration $\sigma \in \{\pm 1\}^N$ such that $H_N(\sigma)$ is at least a fixed positive fraction of the maximal value, no matter how small we allow the fraction to be~\cite{arora2005non}. But here we depart from such worst-case analysis, and focus instead on ``typical'' choices of the coefficients $(W_{ij})$, by postulating that they are chosen randomly. 

Standard concentration inequalities allow us to show that the maximum in \eqref{e.def.max} deviates only little from its expectation in the limit of large $N$, so we may as well focus on studying its expectation. In addition to the expectation of \eqref{e.def.max}, it is natural to also consider, for each $\beta \ge 0$, the quantity
\begin{equation}  
\label{e.def.FN}
F_N(\beta) := \frac 1 N \E \log \bigg(\frac 1 {2^N}\sum_{\si \in \{-1,1\}^N} \exp (\beta H_N(\sigma))\bigg).
\end{equation}
There are several reasons for this. 
From the point of view of statistical mechanics, this quantity is closely related to the Gibbs measure associated with $H_N(\sigma)$, which is the probability measure that attributes a probability proportional to $\exp(\beta H_N(\sigma))$ to each configuration $\sigma \in \{-1,1\}^N$. (See \cite[Section 1.1]{HJbook} for some motivations behind the concept of Gibbs measures.) Mathematically, the free energy can also be seen as a sort of Laplace transform of the quantity of interest, and contains much information about the geometry of $H_N$ and the structure of the Gibbs measure. It is often more convenient to work with, and if we so wish, we can a posteriori recover information about the maximum in \eqref{e.def.max} by considering $F_N(\beta)/\beta$ for large $\beta$. The normalization of $H_N$ has been chosen so that $\max_{\{-1,1\}^N} H_N$ is of order $N$. The free energy thus allows us to interpolate between a large-$\beta$ regime in which the sum in \eqref{e.def.FN} is dominated by the contribution of the configurations $\sigma$ for which $H_N(\sigma)$ is large (``energy dominates''), and a small-$\beta$ regime in which the very large number of configurations with relatively small $H_N(\sigma)$ provide the dominant contribution (``entropy dominates'').

The problem of identifying the large-$N$ limit of the free energy $F_N(\beta)$ turns out to be surprisingly rich and difficult. An initial guess for this limit was proposed in the original paper \cite{sherrington1975solvable} that introduced the model, but it was already understood there that the proposed answer could not be valid for large values of $\beta$. In 1979, Giorgio Parisi then came up with a sophisticated non-rigorous procedure, called the replica method, that led to what is now called the \emph{Parisi formula} for this limit \cite{parisi1979infinite, parisi1980order, parisi1980sequence, parisi1983order} (see also \cite{mourrat2024informal, parisi2023nobel} for more on the replica method). After many years of effort, Francesco Guerra and Michel Talagrand managed to prove the Parisi formula rigorously \cite{gue03, Tpaper} in 2003. A more conceptual proof, centered around the fact that the associated Gibbs measure is asymptotically ultrametric\footnote{Although the concept of ultrametricity is central to the topic, it will not be discussed much in this note. The interested reader can consult \cite{mourrat2024informal} for a light introduction, and \cite[Section~5.7]{HJbook} and \cite{pan} for more precision.}, was then developed by Dmitry Panchenko, and generalized to a broader class of models \cite{pan.aom, pan}. 
The Parisi formula takes the following form. 

\begin{theorem}[Parisi formula \cite{gue03, pan.aom, pan, Tpaper}]
\label{t.parisi}
For every $\beta \ge 0$, we have
\begin{equation}
\label{e.parisi}
f(\beta) := \lim_{N\to+\infty}F_N(\beta)=\inf_{\mu \in \mcl P([0,1])}\Ll( \Phi_\mu(0,0) -{\beta^2}\int_0^1 t\mu([0,t])\d t\Rr),
\end{equation}
where $\mcl P([0,1])$ denotes the space of probability measures on $[0,1]$, and $\Phi_\mu :[0,1]\times \R \to \R$ is the solution to 
\begin{equation}\label{e.parisi.pde.dis}
\begin{cases}
-\partial_t\Phi_\mu(t,x)={\beta^2}\Big(\partial_x^2\Phi_\mu(t,x)+\mu([0,t])\big(\partial_x\Phi_\mu(t,x)\big)^2\Big) & \text{for } (t,x) \in [0,1]\times \R,\\
\Phi_\mu(1,x)=\log\cosh(x)& \text{for } x\in \R.
\end{cases}
\end{equation}
\end{theorem}

The formula \eqref{e.parisi} came as a surprise, and its validity was initially controversial. Compared to more classical problems of statistical mechanics, several aspects stand out. The first is simply the complexity of the formula, as the optimization variable is a probability measure, while for more classical models it usually ranges in a finite-dimensional space. A second and perhaps more fundamental surprise is that the limit free energy is expressed as an infimum, rather than a supremum. Indeed, even before passing to the limit $N \to +\infty$, the free energy of essentially any system can be written as a supremum of a functional involving intuitive energy and entropy terms. To be precise, for every probability measure~$\mu$ over a measure space $E$ and every bounded measurable function $g : E \to \R$, we have
\begin{equation}  
\label{e.gibbs.var}
\log \int e^g \, \d \mu = \sup_{\nu \in \mcl P(E)} \Ll( \int g \, \d \nu - \msf H(\nu \, | \, \mu) \Rr), 
\end{equation}
where $\msf H(\nu \, | \, \mu)$ stands for the relative entropy of $\nu$ with respect to $\mu$, which is $+\infty$ if $\nu$ is not absolutely continuous with respect to $\mu$, and is otherwise given by
\begin{equation*}  
\msf H(\nu \, | \, \mu) := \int \frac{\d \nu}{\d \mu} \log \Ll( \frac{\d \nu}{\d \mu} \Rr) \, \d \mu. 
\end{equation*}
The optimum in \eqref{e.gibbs.var} is achieved by the corresponding Gibbs measure, that is, the probability measure $\nu^*$ whose Radon-Nikodym derivative with respect to $\mu$ is proportional to $e^g$. The first term on the right side of \eqref{e.gibbs.var} expresses the average of $g$ under the measure $\nu^*$, while the second term expresses the cost for samples from $\mu$ to look like samples from $\nu^*$. 
In our case, we think of $g$ as being $\beta H_N$ and of $\mu$ as being the uniform measure on $\{-1,1\}^N$. For finite $N$, the representation coming from \eqref{e.gibbs.var} is complicated since $g$ is random and we need to maximize over all probability measures over $\{-1,1\}^N$, but one could a priori hope that some simplifications occur in the limit of large $N$, so that we ultimately end up with a simple formula for the limit free energy of the form
\begin{equation}  
\label{e.fbeta.rep}
f(\beta) = \sup_{(e,s) \in I} (\beta e - s),
\end{equation}
for some explicit set $I \subset \R^2$. In simpler systems of statistical mechanics, one can often obtain a representation of the form in \eqref{e.fbeta.rep} by identifying all pairs $(e,s)$ such that, roughly speaking,
\begin{equation}  
\label{e.energy-entropy}
\frac 1  N \log \bigg(2^{-N} \sum_{\sigma \in \{-1,1\}^N} \1_{\{ H_N(\sigma) \simeq N e \}}\bigg) \simeq -s \, ;
\end{equation}
or at a minimum, one can deduce legible information in the spirit of \eqref{e.energy-entropy} from the identification of the limit free energy.
For spin glasses however, although we know that a representation of the form \eqref{e.fbeta.rep} exists simply because $f$ is convex, I am not aware of a reasonably direct and concrete way to describe the set $I$, in other words the set of energy-entropy pairs that are achievable by the system. This is related to the fact that, even though the convexity of $f$ is easily checked as $F_N$ itself is convex (e.g.\ by \eqref{e.gibbs.var}), it is not at all clear to verify this convexity property directly from the limit expression given by Theorem~\ref{t.parisi}. 

Since maximization problems are much more standard representations of free energies, one may call the variational problem in \eqref{e.parisi} an ``inverted'' variational representation\footnote{The real benefit of this terminology is that it is robust to the sign convention we choose, as many authors (and we too later in this note) would add a minus sign to the definition of the free energy.}.

%
%
%
%
%
%
\section{Un-inverting the Parisi formula}
\label{s.uninverting}

An ``un-inverted'' representation of the limit free energy of the SK model, that is, one that takes the form of a supremum, was recently found. Recall that we denote by $f(\beta)$ the large-$N$ limit of the free energy $F_N(\beta)$ defined in \eqref{e.def.FN}. 

\begin{theorem}[Un-inverted Parisi formula \cite{mourrat2024uninverting}]  
\label{t.uninverted}
Let $(B_t)_{t \ge 0}$ denote a Brownian motion defined on some filtered probability space $(\Omega, \mcl F, (\mcl F_t)_{t \ge 0}, \PP)$, and let $\mart$ denote the space of bounded martingales on $\Omega$. For every $\beta \ge 0$, we have 
\begin{equation}  
\label{e.uninverted}
f(\beta) = \sup_{\al \in \mart} \bigg\{ \beta\sqrt{2}\EE [\alpha_1 B_1]- \EE[\phi^*(\alpha_1)]  - {\beta^2} \sup_{t \in [0,1]} \int_t^1 (s - \EE[\alpha_s^2]) \, \d s \bigg\} ,
\end{equation}
where for every $\lambda \in \R$, we set
\begin{equation*}  
\label{e.def.phi*}
\phi^*(\lambda) = 
\Ll|
\begin{array}{ll}  
\frac 1 2 \Ll[(1+\lambda)\log(1+\lambda) + (1-\lambda)\log(1-\lambda)\Rr] & \text{ if } |\lambda| \le 1, \\
+\infty & \text{ otherwise}.
\end{array}
\Rr.
\end{equation*}
\end{theorem}
The representation in \eqref{e.uninverted} was obtained by manipulating the Parisi formula from Theorem~\ref{t.parisi}, using the fact from \cite{auffinger2015parisi} that the functional inside the infimum in \eqref{e.parisi} is convex, together with duality arguments. The connection between this representation and the finite-$N$ system remains to be discovered. There are however some indications for how this connection might emerge. For instance, the quantity $\EE[\phi^*(\alpha_1)]$ resembles a relative entropy as in \eqref{e.gibbs.var}. Indeed, denoting by
\begin{equation*}  
\msf{Ber}(m) := \frac{1+m}{2}\de_1 + \frac{1-m}{2}\de_{-1}
\end{equation*}
the law of a random variable taking values in $\{-1,1\}$ with mean $m$,
we have that 
\begin{equation*}  
\frac 1 N  \msf H \Ll( \bigotimes_{i = 1}^N \msf{Ber}(m_i) \, | \, (\msf{Ber}(0))^{\otimes N} \Rr)  = \frac 1 N \sum_{i = 1}^N \phi^*(m_i),
\end{equation*}
and this quantity can be rewritten as $\EE[\phi^*(\alpha_1)]$ provided that the law of $\alpha_1$ is $\frac 1 N \sum_{i = 1}^N \delta_{m_i}$. We recall that for our model, we think of the formula \eqref{e.gibbs.var} with $\mu$ chosen to be the uniform measure over $\{-1,1\}^N$, which is $(\msf{Ber}(0))^{\otimes N}$.

In truth, the Gibbs measure (which is the optimizer of \eqref{e.gibbs.var}) is not a product measure of the form $\bigotimes_{i = 1}^N \msf{Ber}(m_i)$, and the last term in \eqref{e.uninverted} accounts for a correction. For specialists, the structure of the formula in \eqref{t.uninverted} is likely to evoke the finite-$N$ representation of the free energy first introduced by Thouless, Anderson and Palmer \cite{thouless1977solution} and further developed in \cite{auffinger2019spin, auffinger2019thouless, chen2018tap, chen2023generalized}, where a similar correction term appears.

To give further substance to the claim that there should be a direct way to understand the emergence of the variational formula in \eqref{e.uninverted} from the finite-$N$ system, I would like to make a detour to questions that concern optimization algorithms; a more detailed overview of these developments is in \cite{gamarnik2025turing}. First, by taking the large-$\beta$ limit of $f(\beta)/\beta$, one can show that 
\begin{multline}  
\label{e.uninverted.max}
\lim_{ N \to +\infty} \frac 1 N \E \max_{\sigma \in \{-1,1\}^N} H_N(\sigma) 
\\
= \sup_{\alpha \in \mart} \Ll\{\sqrt{2}\EE[\alpha_1 B_1] \ : \ |\alpha_1| \le 1 \text{ and } \forall t \in [0,1], \ \int_t^1 (\EE[\alpha_s^2] -s) \, \d s \ge 0 \Rr\}.
\end{multline}
In a series of recent works, an algorithmic threshold $\msf{ALG}$ was identified such that the following holds. On the one hand, there exists an efficient algorithm that, given the coefficients $(W_{ij})$,  returns a configuration $\si$ such that $H_N(\sigma)/N$ is $\msf {ALG}(1+o(1))$ with probability tending to $1$ as $N$ tends to infinity \cite{elalaoui2021optimization, huang2023algorithmic, jekel2025potential, montanari2021optimization, sellke2021optimizing, subag2021following}. On the other hand, no matter how small $\eps > 0$ is chosen, with probability tending to $1$ as $N$ tends to infinity, an algorithm that is Lipschitz continuous with respect to the input weights $(W_{ij})$ will not be able to output a configuration $\sigma$ such that $H_N(\sigma)/N$ exceeds $\msf{ALG} - \eps$ \cite{gamarnik2021overlap-survey, gamarnik2021overlap-paper, gamarnik2022disordered, huang2025tight}. This threshold value $\msf {ALG}$ can be written as
\begin{equation}  
\label{e.alg.formula}
\msf {ALG} = \sup_{\alpha \in \mart} \Ll\{\sqrt{2}\EE[\alpha_1 B_1] \ : \ |\alpha_1| \le 1 \text{ and } \forall t \in [0,1], \  \EE[\alpha_t^2] = t  \Rr\}.
\end{equation}
One point I find very interesting is that in this case, we have a clear connection between the variational formula \eqref{e.alg.formula} and the spin-glass system at finite $N$. Indeed, for essentially every choice of martingale $\alpha \in \mart$ that satisfies the constraints $|\alpha_1| \le 1$ and $\forall t,\, \EE[\alpha_t^2] = t$, one can construct an algorithm that outputs a configuration $\sigma$ such that $H_N(\sigma)/N$ is approximately $\sqrt{2}\EE[\alpha_1 B_1]$. The algorithm iteratively updates a point in $\R^N$, with small increments that we can here approximate by a continuous evolution $m : [0,1] \to \R^N$, and for each $t \in [0,1]$, the empirical measure of the coordinates $\frac 1 N \sum_{i = 1}^N \delta_{m_i(t)}$ converges weakly to the law of $\alpha_t$. 

In short, I feel particularly interested in, and have the impression that we can make progress upon, the following question:

Can we interpret the variational formula in \eqref{e.uninverted} in terms of finite-$N$ constructs?

Part of my interest in this question is that I think that it has the potential to give us a new way to understand and study spin glasses. Also, once we have a satisfactory answer to this question at a heuristic level, perhaps we could devise a new proof at least of the inequality stating that, for every $\alpha \in \mart$,
\begin{equation}  
\label{e.desired.ineq}
f(\beta) \ge \beta\sqrt{2}\EE [\alpha_1 B_1]- \EE[\phi^*(\alpha_1)]  - {\beta^2} \sup_{t \in [0,1]} \int_t^1 (s - \EE[\alpha_s^2]) \, \d s.
\end{equation}
As will be explained, having a direct proof of the inequality \eqref{e.desired.ineq} that does not rely on the Parisi formula would be tremendously useful for resolving the main open problem discussed in the next section.

%
%
%
%
%
%
\section{Towards more general models}
\label{s.general}

Building on the insights provided by the Parisi formula and its proof, along with subsequent developments, we now have a much deeper understanding of the Sherrington-Kirkpatrick model and the structure of its Gibbs measures.
Moreover, the ideas first developed for the SK model have proven remarkably fruitful in a wide range of other contexts that exhibit similar kinds of ``frustration'', ranging from statistics and high-dimensional geometry to computer science and combinatorics.
Examples include random constraint satisfaction problems \cite{ding2015proof, ding2016satisfiability,  krzakala2007gibbs, mezmon, mezard2002analytic, monasson1999determining}, 
community detection and related large-scale statistical learning  problems \cite{abbe2017community, zdekrz} \cite[Section~4]{HJbook}, error-correcting codes in information theory \cite{richardson2008modern}, and classical combinatorial problems such as graph coloring \cite{coja2018, ding2016maximum, mulet2002coloring}.

\begin{wrapfigure}{r}{0.3\textwidth}
\centering

\includegraphics[width=0.25\textwidth]{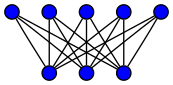}
\captionsetup{width=0.95\linewidth}
\caption{{\small The graph of direct interactions in the bipartite model.}}
    \label{f.bipartite}
\end{wrapfigure}

This being said, and perhaps surprisingly, some models that seem like modest generalizations of the SK model still resist analysis. One such example can be constructed as follows. In the definition of $H_N$ for the SK model in \eqref{e.def.HN}, we sum over all pairs $(i,j)$. For the new model we consider now, we imagine that the indices are organized over two layers as on Figure~\ref{f.bipartite}, and we only sum over pairs of indices that belong to different layers. Such models are related to several classical models of artificial neural networks, such as the Hopfield model \cite{amari1972learning, amit1985spin, amit1987statistical, barra2012equivalence, hopfield1982neural, little1974existence, Tbook2} and restricted Boltzmann machines \cite{barra2018phase, hinton2012practical, hinton2025nobel, smolensky1986information, tubiana2018restricted, tubiana2017emergence}.

To formalize the model precisely, we can first give ourselves, for each integer $N$, two integers $N_1(N)$ and $N_2(N)$ that represent the sizes of the two layers displayed on Figure~\ref{f.bipartite}. We will keep the dependency of $N_1$ and $N_2$ on $N$ implicit from now on, and assume that there exist $\lambda_1, \lambda_2 \in (0,+\infty)$ such that 
\begin{equation}
\label{e.def.lambda}
\lim_{N \to \infty} \frac{N_1}{N} = \lambda_1 \quad \text{ and } \quad \lim_{N \to \infty} \frac{N_2}{N} = \lambda_2. 
\end{equation}
Now, for each $\sigma = (\sigma_1, \sigma_2) = (\sigma_{1,1},\ldots, \sigma_{1,N_1}, \sigma_{2,1}, \ldots, \sigma_{2,N_2}) \in \R^{N_1} \times \R^{N_2}$, we set
\begin{equation}  
\label{e.def.HN.bip}
H_N^\mathrm{bip}(\sigma) := \frac 1 {\sqrt{N}} \sum_{i\le N_1, j \le N_2} W_{i,j} \sigma_{1,i} \sigma_{2,j}.
\end{equation}
We refer to this model as the \emph{bipartite model}.  
Writing $\Sigma_N := \{-1,1\}^{N_1} \times \{-1,1\}^{N_2}$, we would like for instance to understand the large-$N$ behavior of
\begin{equation}  
\label{e.bip.free}
\frac 1 N \E \log \bigg(\frac 1 {|\Sigma_N|}\sum_{\sigma \in \Sigma_N} \exp(\beta H_N^\mathrm{bip}(\sigma))\bigg).
\end{equation}

While this bipartite model may seem like a small modification of the SK model, to this day, we do not know what the limit of the quantity in \eqref{e.bip.free} is; in fact, we do not even know that it converges as $N$ tends to infinity in this case. (The same goes for the maximum of $H_N^\bip$ over $\Sigma_N$.) 
Moreover, the problem we encounter here goes beyond that of adjusting some technical part of the proof of the Parisi formula. Indeed, although one can a priori imagine several possible ways of extending the Parisi formula to this bipartite model, one can in fact show that none of those candidates for the limit are valid \cite[Section~6]{mourrat2020nonconvex}. 

In order to clarify the key difference between the bipartite and the SK models at the technical level, it is useful to change a bit our viewpoint on the definition of these random fields $H_N$ and $H_N^\mathrm{bip}$. Instead of writing them down explicitly as in \eqref{e.def.HN} and \eqref{e.def.HN.bip}, we can equivalently specify that they are centered Gaussian fields, and display their covariance. For the SK model, we have for every $\sigma, \tau \in \R^N$ that
\begin{equation}  
\label{e.cov.HN}
\E \Ll[ H_N(\sigma) H_N(\tau) \Rr] = N \Ll( \frac{\sigma \cdot \tau}{N} \Rr) ^2,
\end{equation}
where $\sigma \cdot \tau$ denotes the scalar product between $\sigma$ and $\tau$. More generally, one could consider centered Gaussian fields $(H_N(\sigma))_{\sigma \in \R^N}$ such that, for some smooth function $\xi : \R \to \R$, we have for every $\sigma, \tau \in \R^N$ that
\begin{equation}  
\label{e.def.cov}
\E \Ll[ H_N(\sigma) H_N(\tau) \Rr] = N \xi\Ll( \frac{\sigma \cdot \tau}{N} \Rr);
\end{equation}
this corresponds to an assumption on the invariance of the law of $H_N$ under orthogonal transformations.
The SK model \eqref{e.def.HN} corresponds to the case when $\xi(r) = r^2$. For $\xi(r) = r^3$, we can construct a Gaussian field that satisfies \eqref{e.def.cov} by setting
\begin{equation*}  
H_N(\sigma) := \frac 1 N \sum_{i,j,k = 1}^N W_{i,j,k} \sigma_i \sigma_j \sigma_k,
\end{equation*}
where $(W_{i,j,k})$ are independent centered Gaussians with unit variance. For every integer $p \ge 1$, we can generalize this and construct a centered Gaussian field $H_N$ such that \eqref{e.def.cov} holds with $\xi(r) = r^p$. By considering linear combinations of independent versions of such fields, we can build a centered Gaussian field $H_N$ such that \eqref{e.def.cov} holds with
\begin{equation}  
\label{e.valid.xi}
\xi(r) = \sum_{p =0}^{+\infty} a_p^2 \,  r^p,
\end{equation}
provided that the sequence $(a_p)_{p \ge 1}$ decays to zero sufficiently fast. It turns out that the functions of the form in \eqref{e.valid.xi} are all those such that \eqref{e.def.cov} holds for some centered Gaussian field $H_N$ (see \cite{schoenberg1942positive}, and \cite[Proposition~6.6]{mourrat2023free} for a more general statement covering cases with multiple types of spins).

In the case of the bipartite model \eqref{e.def.HN.bip}, we have instead that, for every $\sigma, \tau \in \R^{N_1} \times \R^{N_2}$,
\begin{equation}  
\label{e.def.cov.bip}
\E \Ll[ H_N^\bip(\sigma) H_N^\bip(\tau) \Rr] = N \Ll( \frac{\sigma_1 \cdot \tau_1}{N} \Rr) \Ll( \frac{\sigma_2 \cdot \tau_2}{N} \Rr).
\end{equation}
The key technical difference between the SK and the bipartite models is that here the relevant function that shows up on the right side of \eqref{e.def.cov.bip} is the mapping $(x,y) \mapsto xy$, which is \emph{not convex}. To be precise, for models with only one type of spins, i.e.\ of the form in \eqref{e.def.cov}, what is crucial is that the function $\xi$ is convex over $\R_+$; as one can see from~\eqref{e.valid.xi}, this is in fact always the case! This convexity property can however break down as soon as we consider models with two or more types of spins. In general, we can consider models with a fixed number $D$ of types of spins, say $\sigma = (\sigma_1, \ldots, \sigma_D) \in \R^{N_1} \times \cdots \times \R^{N_D}$, with $N_d / N \to \lambda_d \in (0,+\infty)$ for every $d \in \{1,\ldots, D\}$, and with a covariance such that, for every $\sigma, \tau \in \R^{N_1} \times \cdots \times \R^{N_D}$,
\begin{equation}  
\label{e.cov.general}
\E \Ll[ H_N(\sigma) H_N(\tau) \Rr] = N \xi \Ll( \Ll( \frac{\sigma_d \cdot \tau_{d}}{N} \Rr)_{1 \le d \le D}  \Rr) ,
\end{equation}
where $\xi$ is some (admissible) function from $\R^{D}$ to $\R$. Those models for which we can write down and rigorously prove a Parisi formula for the limit free energy are those for which the function $\xi$ is convex over $\R^D_+$ \cite{barcon, chen2024on, chen2024free, pan.aom, pan.multi, pan.potts, pan.vec}. Some particular models of the form \eqref{e.cov.general} with $\xi$ that is not convex over $\R^D_+$ but with several additional symmetries have also been successfully analyzed \cite{aufchebi, baik2020free, bates2025balanced, dartois2024injective, subag2023tap2, subag2025tap1}, but for the most part, the analysis of models of the form \eqref{e.cov.general} with non-convex $\xi$ remains open.

%
%
%
%
%
%
\section{A connection with Hamilton-Jacobi equations}
\label{s.pdes}

In order to make progress on the identification of the limit free energy for systems such as the bipartite model, my collaborators and I have explored an approach that consists in seeking a partial differential equation that would be solved by the limit of $F_N$. The free energy $F_N$ as we defined it here in \eqref{e.def.FN} (or in \eqref{e.bip.free} for the bipartite model) depends only on $\beta$, and it is not possible to find a simple equation for $F_N$ or its limit that would only involve derivatives in $\beta$. Hence, we first seek to add terms to the energy function~$H_N$ that depend on additional parameters; for instance, we could replace $\beta H_N(\sigma)$ by $\beta H_N(\sigma) + \lambda H_N'(\sigma)$ for some free parameter $\lambda \in \R$ and some well-chosen $H_N'$. This would yield a free energy that now depends on $\lambda$ in addition to $\beta$. We perform this ``enrichment'' of the free energy with the hope of finding a partial differential equation involving derivatives in, say, $\beta$ and $\lambda$, that the free energy would asymptotically solve as $N$ tends to infinity. Naturally, one would like the additional quantities such as $H_N'$ in the example to be less complicated to analyze than the original field $H_N$. On the other hand, we would like the additional parameters to be sufficiently rich that we can ultimately close the equation for the limit free energy. In practice, we will always shoot for first-order partial differential equations, so we can intuitively think of the task as that of building a simpler but ``locally equivalent'' energy function $H_N'$, so that we can compensate small variations of $\beta$ with small variations of $\lambda$ and keep the free energy roughly constant. If this is indeed possible, then we obtain a way to flow the parameter $\beta$ from the ``easy'' case with $\beta = 0$ towards the value of $\beta$ of interest.

The idea of thinking of the limit free energy of a model of statistical mechanics as a solution to a partial differential equation goes back at least to \cite{bra83, new86} (see also \cite{bauerschmidt2023stochastic} for a recent survey on related topics). For simpler models of statistical mechanics, as well as some problems of inference such as community detection on dense graphs, this strategy can work very well, see \cite{chen2022statistical} and \cite[Chapters 1 to 4]{HJbook} for a detailed presentation. The case of spin glasses is more difficult though. Ideas in this spirit were first explored in \cite{abarra, barra2,barra1, guerra2001sum} under simplifying assumptions. We now informally discuss some of the recent progress in this direction. Similar difficulties, not discussed further here, also show up for community detection on sparse graphs \cite{sparse_PDE, sparse_prob, mutual_information}. 

In order to keep the notation simple, we first present the approach in the case of the SK model \eqref{e.def.HN}. For every $t \ge 0$ and $h \ge 0$, we set
\begin{equation}  
\label{e.def.newFN}
F_N(t,h) := - \frac 1 N \E \log \bigg(\frac 1 {2^N} \sum_{\sigma \in \{- 1,1\}^N} \exp(\sqrt{2t} H_N(\sigma) - Nt + \sqrt{2h} z \cdot \sigma - N h)\bigg), 
\end{equation}
where $z = (z_1,\ldots, z_N)$ is a vector of independent centered Gaussian random variables with unit variance, independent of $H_N$, and we recall that the function $H_N$ for the SK model is defined in~\eqref{e.def.HN}. 
The key decision we made is that of adding the term involving $z \cdot \sigma$. As will be seen more clearly below, the fact that this term is linear in $\sigma$ makes it much simpler indeed. And if we write $H_N(\sigma)$ in the form 
\begin{equation*}  
H_N(\sigma) = \frac 1 {\sqrt N} \sum_{i = 1}^N \Ll( \sum_{j = 1}^N W_{ij} \sigma_j\Rr) \sigma_i ,
\end{equation*}
it is perhaps not unreasonable to hope that the term $(\sum_{j = 1}^N W_{ij} \sigma_j)$ could be substituted with equivalent independent Gaussians. We also wrote a factor of $\sqrt{2t}$ in place of $\beta$ in front of $H_N$; since $H_N$ is Gaussian, this ensures that the variance of $\sqrt{2t} H_N$ scales linearly, as with Brownian motion. The compensating parameter $Nt$ is only a convenience\footnote{In general, we do not need to restrict ourselves to models defined on $\{-1,1\}^N$. For a model whose covariance is given by \eqref{e.def.cov}, we can consider
\begin{equation*}  
F_N(t,h) := - \frac 1 N \E \log \int \exp\Ll(\sqrt{2t} H_N(\sigma) - N t\xi(|\sigma|^2/N) + \sqrt{2h} z \cdot \sigma - h|\sigma|^2\Rr) \, \d P_N(\sigma),
\end{equation*}
where $P_N = P_1^{\otimes N}$ is the $N$-fold tensor product of a probability measure $P_1$ on $\R$ with compact support. We recover the SK model \eqref{e.def.HN} by choosing $\xi(r) = r^2$ and $P_1 = (\de_1 + \de_{-1})/2$. Notice that with this definition (and thanks to the minus sign appearing there), Jensen's inequality yields that $F_N \ge 0$. The additional terms involving $t\xi(|\sigma|^2/N)$ and $|\sigma|^2$ facilitate the analysis, and can be removed a posteriori as they are not themselves random. 
}.
Similar comments also hold concerning the terms $\sqrt{2h}$ and $-Nh$. Notice also that, unlike in previous sections, we added a minus sign to the definition of the free energy $F_N(t,h)$.

Before proceeding, we introduce notation for the Gibbs measure. For any function $f$, we write
\begin{equation}  
\label{e.def.Gibbs}
\langle f(\sigma) \rangle := \frac{\sum_{\sigma \in \{\pm 1\}^N} f(\sigma) \exp(H_N(t,h,\sigma))}{\sum_{\sigma \in \{\pm 1\}^N} \exp(H_N(t,h,\sigma))},
\end{equation}
where  $H_N(t,h,\sigma) := \sqrt{2t} H_N(\sigma) - Nt + \sqrt{2h} z \cdot \sigma - N h$. In the notation on the left side of \eqref{e.def.Gibbs}, the bracket $\langle \cdot \rangle$ stands for the expectation with respect to the Gibbs measure, and we think of $\sigma$ as a random variable that is sampled accordingly. We write $\sigma'$ to denote an independent copy of $\sigma$ under the Gibbs measure, so that
\begin{equation*}  
\langle f(\sigma, \sigma') \rangle := \frac{\sum_{\sigma,\sigma' \in \{\pm 1\}^N} f(\sigma, \sigma') \exp(H_N(t,h,\sigma)+H_N(t,h,\sigma'))}{\sum_{\sigma,\sigma' \in \{\pm 1\}^N} \exp(H_N(t,h,\sigma)+H_N(t,h,\sigma'))}.
\end{equation*}
This expectation $\langle \cdot \rangle$ depends on the parameters $t$ and $h$, even though we keep it implicit in the notation. 
A simple calculation involving Gaussian integration by parts gives us that
\begin{equation}  
\label{e.FN.derivatives}
\dr_t F_N(t,h) = \E \la \Ll( \frac{\sigma \cdot \sigma'}{N} \Rr) ^2 \ra \quad \text{ and } \quad \dr_h F_N(t,h) = \E \la  \frac{\sigma \cdot \sigma'}{N}  \ra  . 
\end{equation}
For general models as in \eqref{e.def.cov}, we would find the same expression for $\dr_h F_N$ as in \eqref{e.FN.derivatives} (for the corresponding definition of the Gibbs measure), while for the derivative in $t$, we would find that
\begin{equation}  
\label{e.drt.xi}
\dr_t F_N(t,h) = \E \la \xi \Ll( \frac{\sigma \cdot \sigma'}{N} \Rr)  \ra.
\end{equation}
Coming back to the SK model for now, we thus obtain that 
\begin{equation}
\label{e.pde.FN}
\dr_t F_N - (\dr_h F_N)^2 = \E \la \Ll( \frac{\sigma \cdot \sigma'}{N} \Rr) ^2 \ra  - \Ll( \E \la  \frac{\sigma \cdot \sigma'}{N}  \ra   \Rr) ^2.
\end{equation}
The right-hand side of \eqref{e.pde.FN} is the variance of the random variable $\sigma \cdot \sigma'/{N}$ under $\E \la \cdot \ra$. Since $\sigma \cdot \sigma'/N$ is a sum of a large number of terms, we may at first anticipate that it will have small fluctuations. If we assume that this is so for the moment, we are led to the expectation that $F_N$ may converge to a limit function~$f$ that solves the equation
\begin{equation}  
\label{e.hj.SK.simple}
\dr_t f - (\dr_h f)^2 = 0.
\end{equation}
Moreover, one can easily compute the value of the free energy $F_N$ \eqref{e.def.newFN} at $t = 0$, as we have that, for every $h \ge 0$,
\begin{equation}
\label{e.simple.init}
F_N(0,h) = F_1(0,h). 
\end{equation}
Hence, if we believe that the random variable $\sigma \cdot \sigma'/N$ has vanishingly small fluctuations in the limit of large $N$, then we are led to the belief that $F_N$ should converge to the function $f$ that solves \eqref{e.hj.SK.simple} with initial condition $f(0,\cdot) = F_1(0,\cdot)$. For the model with the covariance as in~\eqref{e.def.cov}, we can proceed in the same way, and under the same assumption, we would obtain the limit partial differential equation
\begin{equation}
\label{e.hj.xi.simple}
\dr_t f  - \xi(\dr_h f) = 0.
\end{equation}
Equations of this form go by the name of Hamilton-Jacobi equations. 

Unfortunately, the hypothesis that the random variable $\sigma \cdot \sigma'/N$ has vanishingly small fluctuations in the limit of large $N$ is only valid at high temperature, or in other words, for small values of $t$. For large values of $t$, the Gibbs measure becomes more complex (with an ultrametric structure), and the variance of $\sigma \cdot \sigma'/N$ under the Gibbs measure does not tend to zero. 

Since for large $t$, one cannot close the equation for the limit of $F_N$ using only the variables $t$ and $h$, we need to refine this first attempt and introduce a richer additional term than this $z \cdot \sigma$ that we used here. The more sophisticated term is still linear in $\sigma$, but replaces the simple ``external field'' $z$ by one with an ultrametric structure. This ultrametric structure is encoded by a number of parameters, which can be collectively bundled into an increasing\footnote{By ``increasing'', we mean that for every $s \le t \in [0,1)$, we have $q(s) \le q(t)$.} and bounded cadlag function $q : [0,1) \to \R$; we denote by $\mcl Q$ the set of such functions. 
The detailed motivation and complete definition of this term would bring us too far off, and I will have to ask the reader to accept (or to consult \cite[Section~6]{HJbook}) that it is indeed possible to define an enriched free energy $F_N(t, q)$, for every $t \in \R_+$ and $q \in \mcl Q$, so that asymptotically as $N$ tends to infinity, we have for the SK model that
\begin{equation}  
\label{e.approx.hj}
\partial_t F_N - \int_0^1 (\partial_q F_N)^2  = \text{some plausibly small conditional variance term},
\end{equation}
and that moreover, we have in analogy with \eqref{e.simple.init} that, for every $q \in \mcl Q$,
\begin{equation}  
\label{e.def.psi}
F_N(0,q) = F_1(0,q). \quad \text{ For convenience, we write } \psi_1(q) := F_1(0,q). 
\end{equation}
In the special case when the path $q$ is identically equal to $h$, we recover the quantity in~\eqref{e.def.newFN}, so this enriched free energy $F_N : \R_+ \times \mcl Q \to \R$ is an extension of that defined in \eqref{e.def.newFN}; and in particular, the quantity we are ultimately most interested in computing is $F_N(t,0)$. Informally, the derivative $\partial_q$ appearing in  \eqref{e.approx.hj} is such that, for a sufficiently smooth function $g : \mcl Q \to \R$ and $q \in \mcl Q$, the quantity $\partial_q g(q,\cdot)$ is a function from $[0,1]$ to $\R$ such that for every $q' \in \mcl Q$ and as $\eps$ tends to zero, 
\begin{equation*}  
g((1-\ep) q + \ep q') - g(q) = \ep \int_0^1 \partial_q g(q,u) (q'-q)(u) \, \d u + o(\ep) \, ;
\end{equation*}
and a more explicit writing of the integral in \eqref{e.approx.hj} is 
$
\int_0^1 (\partial_q F_N(t,q,u))^2 \, \d u.
$

For all models with a single type, we can indeed characterize the limit free energy as the unique solution to the equation in \eqref{e.approx.hj}. 
\begin{theorem}[Limit free energy via a Hamilton-Jacobi equation \cite{chen2022hamilton, mourrat2022parisi, mourrat2020extending}]
\label{t.hj}
The enriched free energy $F_N : \R_+ \times \mcl Q \to \R$ of the SK model converges pointwise to the unique function $f : \R_+ \times \mcl Q \to \R$ that solves 
\begin{equation}
\label{e.hj.SK}
\begin{cases}
\dr_t f - \int_0^1 (\dr_q f)^2  = 0   & \text{ on } \R_+ \times \mcl Q, \\
f(0,\cdot) = \psi_1  & \text{ on } \mcl Q.
\end{cases}
\end{equation}
More generally, if $F_N : \R_+ \times \mcl Q \to \R$ stands instead for the enriched free energy associated with a model with covariance given by \eqref{e.def.cov}, then $F_N$ converges pointwise to the unique function $f : \R_+ \times \mcl Q \to \R$ that solves
\begin{equation}
\label{e.hj.xi}
\begin{cases}
\dr_t f - \int_0^1 \xi(\dr_q f)   = 0   & \text{ on } \R_+ \times \mcl Q, \\
f(0,\cdot) = \psi_1  & \text{ on } \mcl Q.
\end{cases}
\end{equation}
\end{theorem}
Part of the task of making sense of this theorem is that one needs to find a good notion of solution for the Hamilton-Jacobi equations in \eqref{e.hj.SK} and \eqref{e.hj.xi}. This is based on the notion of viscosity solutions (see \cite[Chapter 3]{HJbook} for an introduction that is tailored to our context). 

Using that $\xi$ is convex on $\R_+$, we can in fact write the viscosity solution $f$ of \eqref{e.hj.xi} as a variational formula. Indeed, this solution is such that, for every $t \ge 0$ and $q \in \mcl Q$,
\begin{equation}  
\label{e.hopf-lax.xi}
f(t,q) = \sup_{q' \in \mcl Q} \Ll(\psi_1(q+q') - t \int_0^1 \xi^* \Ll( \frac{q'}{t} \Rr)  \Rr),
\end{equation}
where $\xi^*$ is the convex dual of $\xi$, which is defined, for every $s \in \R$, by
\begin{equation*}  
\xi^*(s) := \sup_{r \ge 0} \Ll( rs - \xi(r) \Rr) .
\end{equation*}
One can recover the Parisi formula by setting $q = 0$ in~\eqref{e.hopf-lax.xi}, making a change of variables, and doing some explicit calculations involving the function~$\psi_1$ (see again \cite[Section~6]{HJbook} for more details)\footnote{Recall that in this section, we added a minus sign in the definition of the free energy in \eqref{e.def.newFN}, so what we see as a supremum here is an infimum in the convention of Sections~\ref{s.parisi} and \ref{s.uninverting}.}.

The main motivation for developing this PDE point of view on the Parisi formula is to tackle the case when $\xi$ is non-convex, such as is the case for the bipartite model. In this case, with $H_N^\bip$ defined in~\eqref{e.def.HN.bip}, we can try to mimic the simple arguments that led us to \eqref{e.hj.SK.simple} or \eqref{e.hj.xi.simple}. The point now is that since there are two types of spins, we add one variable for each of these two types. That is, for every $t \ge 0$, $h = (h_1,h_2) \in \R_+^2$, and $\sigma = (\sigma_1, \sigma_2) \in \Sigma_N$, we set
\begin{equation*}  
H_N(t,h,\sigma) := \sqrt{2t} H_N^\bip(\sigma) - N\lambda_1 \lambda_2 t + \sqrt{2h_1} z_1 \cdot \sigma_1 - N_1 h_1 + \sqrt{2 h_2} z_2 \cdot \sigma_2 - N_2 h_2,
\end{equation*}
as well as
\begin{equation*}  
F_N(t,h) := - \frac 1 N \E \log \bigg( \frac{1}{|\Sigma_N|} \sum_{\sigma \in \Sigma_N} \exp(H_N(t,h,\sigma)) \bigg).
\end{equation*}
In the displays above and from now on, we drop the superscript $^\bip$ for ease of notation. 
The definition of the Gibbs average $\la \cdot \ra$ is as in the formula in \eqref{e.def.Gibbs}, with $h \in \R_+^2$ and with the summation variable $\sigma$ ranging in $\Sigma_N$. In place of \eqref{e.FN.derivatives} or \eqref{e.drt.xi}, we obtain that 
\begin{equation*}  
\dr_t F_N(t,h) = \E \la \Ll( \frac{\sigma_1 \cdot \sigma_1'}{N} \Rr) \Ll( \frac{\sigma_2 \cdot \sigma_2'}{N} \Rr) \ra,
\end{equation*}
while we still have
\begin{equation*}  
\dr_{h_1} F_N(t,h) = \E \la  \frac{\sigma_1 \cdot \sigma_1'}{N}  \ra  \quad \text{ and } \quad \dr_{h_2} F_N(t,h) = \E \la  \frac{\sigma_2 \cdot \sigma_2'}{N}  \ra  .
\end{equation*}
Under the assumption that the random variables $\sigma_1 \cdot \sigma_1'/N$ and $\sigma_2 \cdot \sigma_2'/N$ have vanishing fluctuations in the limit of large $N$, we would expect $F_N$ to converge to $f$ solution to
\begin{equation*}  
\dr_t f - \dr_{h_1} f \,  \dr_{h_2} f = 0. 
\end{equation*}
The initial condition $f(0,\cdot)$ is easy to compute as it factorizes similarly to what we saw in \eqref{e.simple.init}. 

As was the case earlier, the assumption of concentration of the random variables $\sigma_1 \cdot \sigma_1'/N$ and $\sigma_2 \cdot \sigma_2'/N$ is invalid for large $t$, and we need again to pass to a more sophisticated free energy $F_N : \R_+ \times \mcl Q^2 \to \R$, which this time takes as second argument a pair of paths $q = (q_1,q_2) \in \mcl Q^2$.
We have a factorization property similar to that in the first identity of \eqref{e.def.psi}; for every $q = (q_1,q_2) \in \mcl Q^2$, we write
\begin{equation*}  
\psi_2(q) := \lim_{N \to +\infty} F_N(0,q).
\end{equation*}
It is useful here to distinguish between the function $\psi_1$ defined on $\mcl Q$ we introduced earlier in \eqref{e.def.psi}, and that new function $\psi_2 : \mcl Q^2 \to \R$ we just defined, as indeed one can easily show using \eqref{e.def.lambda} that
\begin{equation}  
\label{e.psi2.decomp}
\psi_2(q) = \lambda_1 \psi_1(q_1) + \lambda_2 \psi_1(q_2). 
\end{equation}
Calculations similar to those leading to \eqref{e.approx.hj} lead to the following conjecture.
\begin{conjecture}  
\label{conj}
The enriched free energy $F_N : \R_+ \times \mcl Q^2 \to \R$ for the bipartite model converges to the function $f : \R_+ \times \mcl Q^2 \to \R$ that solves
\begin{equation}
\label{e.hj.bip}
\begin{cases}
\dr_t f - \int_0^1 \dr_{q_1} f \, \dr_{q_2} f  = 0   & \text{ on } \R_+ \times \mcl Q^2, \\
f(0,\cdot) = \psi_2  & \text{ on } \mcl Q^2.
\end{cases}
\end{equation}
\end{conjecture}
In equation \eqref{e.hj.bip}, the nonlinearity (i.e. the mapping $(x,y) \mapsto xy$) is neither convex nor concave, so we cannot write a variational representation similar to that in \eqref{e.hopf-lax.xi} for the solution to \eqref{e.hj.bip}.

A number of partial results have been obtained that give credence to Conjecture~\ref{conj}. For ease of discussion, let us suppose that the enriched free energy $F_N : \R_+ \times \mcl Q^2 \to \R$ of the bipartite model converges pointwise to some function $g$ (one can easily show that converging subsequences exist). First, we know from \cite{chen2022hamilton, mourrat2020nonconvex, mourrat2023free} that $g \ge f$, where $f$ is the (unique) viscosity solution to \eqref{e.hj.bip}. Moreover, the limit function $g$ is differentiable ``almost everywhere''\footnote{Quotation marks are due here because there is no Lebesgue measure on $\mcl Q$; the exact formulation is in terms of Gaussian null sets.}, and the function $g$ solves the equation displayed in \eqref{e.hj.bip} at every point of differentiability \cite{chen2024free}. The latter property is weaker than that of being a viscosity solution to \eqref{e.hj.bip} though, so this is not sufficient to conclude. Another result from \cite{chen2024free} gives a characterization of $g$ in terms of a critical point of an explicit functional. While one can show that there exists a unique such critical point for small $t$, there can be more than one critical point for large $t$, so this is again not sufficient to conclude. For those readers who are familiar with Hamilton-Jacobi equations, the result can be understood as saying that the value of $g$ at a given point $(t,q)$ is as prescribed by one of the characteristic lines that goes through $(t,q)$.

%
%
%
%
%
%

\section{Conjectured un-inverted formula for general models}
\label{s.conjecture}

As already mentioned, direct generalizations of the Parisi formula to models such as the bipartite case yield formulas that are provably invalid. For a long time, I therefore could not see any alternative candidate variational formula for the limit of $F_N$. I have changed my mind on this, and now believe that the limit of $F_N$ in the bipartite case does admit a variational representation. However, I do not expect it to take a form similar to that of the Parisi formula in \eqref{e.parisi}, but rather to take a form similar to its ``un-inverted'' version in \eqref{e.uninverted}. 

For Hamilton-Jacobi equations such as \eqref{e.hj.bip}, there are two classical settings in which one can write a variational formula. The first, already discussed, is when the non-linearity in the equation is either convex or concave, but this is clearly not the case here. The second classical situation in which one can obtain a variational formula for the solution to \eqref{e.hj.bip} is when the initial condition $\psi_2$ is either convex or concave. Alas, one can also see that the function $\psi_2$ is neither convex nor concave in general \cite[Section~6]{mourrat2020nonconvex}. 
This being said, the function $\psi_2$ can be turned into a concave function through the following change of variables. Recalling the decomposition in \eqref{e.psi2.decomp}, we first work with the function $\psi_1$, and interpret a single path $q \in \mcl Q$ as the inverse cumulative distribution function of a probability measure. In other words, there is a bijective correspondence between the set $\mcl Q$ and the set $\mcl P_c(\R_+)$ of probability measures over $\R_+$ with compact support, through the mapping which to a given $q \in \mcl Q$ attributes the law of $q(U)$, where $U$ is a uniform random variable over $[0,1]$. Denoting this mapping by $M : \mcl Q \to \mcl P_c(\R_+)$, we define $\td \psi_1 : \mcl P_c(\R_+) \to \R$ so that $\td \psi_1(M q) = \psi_1(q)$. It turns out that the function $\td \psi_1$ is concave over $\mcl P_c(\R_+)$  \cite{auffinger2015parisi}. Similarly, we can define $\td \psi_2 : (\mcl P_c(\R_+))^2 \to \R$ such that for every $q_1,q_2 \in \mcl Q$, we have $\td \psi_2(M q_1, M q_2) = \psi_2(q_1, q_2)$, and by \eqref{e.psi2.decomp}, the function $\td \psi_2$ is concave over $(\mcl P_c(\R_+))^2$. 

Since $\td \psi_2$ is concave, it can be written as an infimum of affine functions. I now expect that the solution to the Hamilton-Jacobi equation in \eqref{e.hj.bip} can be represented as the infimum of the solutions started from these enveloping affine functions. Using also explicit calculations involving the function $\psi_1$, this would yield an ``un-inverted'' representation of the limit free energy in the spirit of the right-hand side of \eqref{e.uninverted}. 

For general models of the form in \eqref{e.cov.general} with convex $\xi$, the connection between the Parisi formula, the Hamilton-Jacobi equation, and the un-inverted variational representation is verified in \cite{chen2025convex}. One natural goal, also interesting from a purely PDE perspective, is to extend the un-inverted variational representation of the Hamilton-Jacobi equation to the non-convex case. Since we know from \cite{chen2022hamilton, mourrat2020nonconvex, mourrat2023free} that the limit free energy is always greater than or equal to the solution to the Hamilton-Jacobi equation, this would give us a lower bound on the free energy in the form of a variational formula similar to the right-hand side of \eqref{e.uninverted}. The justification of an inequality in the spirit of \eqref{e.desired.ineq} would then allow us to complete this picture, prove Conjecture~\ref{conj}, and obtain a variational representation of the limit free energy for all models with a covariance taking the form in~\eqref{e.cov.general}.

%
%
%
%
%
%
\section{Conclusion}
\label{s.concl}

Despite the relative simplicity of their definition, spin glasses display a mathematical structure that I find surprisingly rich and profound. Moreover, the ideas and techniques that were developed to study spin glasses have turned out to be useful in the analysis of a broad class of models across disciplines. 

This short review is centered around the problem of identifying the free energy of models of spin glasses involving several different types of spins. We chose the bipartite model defined in \eqref{e.def.HN.bip} as our guiding example of models in the class of centered Gaussian fields with a covariance in the form of \eqref{e.cov.general}. While the case when $\xi$ is convex over $\R^D_+$ is well-understood, models with non-convex $\xi$ such as the bipartite model have so far resisted complete analysis. Strikingly, direct generalizations of the Parisi formula to these models yield predictions that are demonstrably false, and an alternative approach is required. 

We reviewed one possible approach based on the idea that the limit free energy should satisfy a Hamilton-Jacobi equation. One can formulate the precise Conjecture~\ref{conj} to this effect, and several partial results have been obtained that give substance to it.

Perhaps most interestingly, a connection between the Hamilton-Jacobi equation appearing in Conjecture~\ref{conj} and an ``un-inverted'' variational formula in the spirit of that presented in Theorem~\ref{t.uninverted} is starting to emerge. In my opinion, this ``un-inverted'' formula deserves further study and could potentially help us to better understand spin glasses, including those models for which the Parisi formula is already proved rigorously.

The topic of spin glasses is much broader than what this short and partial review could cover. Books on spin glasses include \cite{bolthausen2007spin, bovier2006statistical, bovier1998mathematical, charbonneau2023spin, contucci2013perspectives, dominicis2006random, HJbook, mpvbook, nishimori2001statistical, opper2001advanced, pan, stein2013spin, Tbook1, Tbook2}. The book \cite{HJbook} is close in spirit to the discussion presented here, and in particular to Section~\ref{s.pdes}.

\bigskip

\noindent \textbf{Acknowledgements.} I would like to warmly thank Hong-Bin Chen, Tom\'as Dominguez, and Victor Issa, with whom much of what is presented here was developed.

\small
\bibliographystyle{plain}
\bibliography{gazette}

\newcommand{\noop}[1]{} \def\cprime{$'$}
\begin{thebibliography}{10}

\bibitem{abbe2017community}
Emmanuel Abbe.
\newblock Community detection and stochastic block models: recent developments.
\newblock {\em J. Mach. Learn. Res.}, 18(1):6446--6531, 2017.

\bibitem{abarra}
Elena Agliari, Adriano Barra, Raffaella Burioni, and Aldo Di~Biasio.
\newblock Notes on the p-spin glass studied via {H}amilton-{J}acobi and
  smooth-cavity techniques.
\newblock {\em J. Math. Phys.}, 53(6):063304, 29, 2012.

\bibitem{amari1972learning}
Shun'ichi Amari.
\newblock Learning patterns and pattern sequences by self-organizing nets of
  threshold elements.
\newblock {\em IEEE Transactions on computers}, 100(11):1197--1206, 1972.

\bibitem{amit1985spin}
Daniel~J. Amit, Hanoch Gutfreund, and Haim Sompolinsky.
\newblock Spin-glass models of neural networks.
\newblock {\em Phys. Rev. A}, 32(2):1007, 1985.

\bibitem{amit1987statistical}
Daniel~J. Amit, Hanoch Gutfreund, and Haim Sompolinsky.
\newblock Statistical mechanics of neural networks near saturation.
\newblock {\em Annals of Physics}, 173(1):30--67, 1987.

\bibitem{arora2005non}
Sanjeev Arora, Eli Berger, Hazan Elad, Guy Kindler, and Muli Safra.
\newblock On non-approximability for quadratic programs.
\newblock In {\em 46th Annual IEEE Symposium on Foundations of Computer Science
  (FOCS'05)}, pages 206--215. IEEE, 2005.

\bibitem{aufchebi}
Antonio Auffinger and Wei-Kuo Chen.
\newblock Free energy and complexity of spherical bipartite models.
\newblock {\em J. Stat. Phys.}, 157(1):40--59, 2014.

\bibitem{auffinger2015parisi}
Antonio Auffinger and Wei-Kuo Chen.
\newblock The {P}arisi formula has a unique minimizer.
\newblock {\em Comm. Math. Phys.}, 335(3):1429--1444, 2015.

\bibitem{auffinger2019spin}
Antonio Auffinger and Aukosh Jagannath.
\newblock On spin distributions for generic {$p$}-spin models.
\newblock {\em J. Stat. Phys.}, 174(2):316--332, 2019.

\bibitem{auffinger2019thouless}
Antonio Auffinger and Aukosh Jagannath.
\newblock Thouless-{A}nderson-{P}almer equations for generic {$p$}-spin
  glasses.
\newblock {\em Ann. Probab.}, 47(4):2230--2256, 2019.

\bibitem{baik2020free}
Jinho Baik and Ji~Oon Lee.
\newblock Free energy of bipartite spherical {S}herrington-{K}irkpatrick model.
\newblock {\em Ann. Inst. Henri Poincar\'e{} Probab. Stat.}, 56(4):2897--2934,
  2020.

\bibitem{barra2012equivalence}
Adriano Barra, Alberto Bernacchia, Enrica Santucci, and Pierluigi Contucci.
\newblock On the equivalence of {Hopfield networks and Boltzmann machines}.
\newblock {\em Neural Networks}, 34:1--9, 2012.

\bibitem{barcon}
Adriano Barra, Pierluigi Contucci, Emanuele Mingione, and Daniele Tantari.
\newblock Multi-species mean field spin glasses. {R}igorous results.
\newblock {\em Ann. Henri Poincar\'{e}}, 16(3):691--708, 2015.

\bibitem{barra2}
Adriano Barra, Gino Del~Ferraro, and Daniele Tantari.
\newblock Mean field spin glasses treated with {PDE} techniques.
\newblock {\em Eur. Phys. J. B}, 86(7):Art. 332, 10, 2013.

\bibitem{barra1}
Adriano Barra, Aldo Di~Biasio, and Francesco Guerra.
\newblock Replica symmetry breaking in mean-field spin glasses through the
  {H}amilton-{J}acobi technique.
\newblock {\em J. Stat. Mech. Theory Exp.}, (9):P09006, 22, 2010.

\bibitem{barra2018phase}
Adriano Barra, Giuseppe Genovese, Peter Sollich, and Daniele Tantari.
\newblock Phase diagram of restricted {B}oltzmann machines and generalized
  {H}opfield networks with arbitrary priors.
\newblock {\em Phys. Rev. E}, 97(2):022310, 2018.

\bibitem{bates2025balanced}
Erik Bates and Youngtak Sohn.
\newblock Balanced multi-species spin glasses.
\newblock {Preprint, arXiv:2507.06522}.

\bibitem{bauerschmidt2023stochastic}
Roland Bauerschmidt, Thierry Bodineau, and Benoit Dagallier.
\newblock Stochastic dynamics and the {P}olchinski equation: an introduction.
\newblock {\em Probab. Surv.}, 21:200--290, 2024.

\bibitem{bolthausen2007spin}
Erwin Bolthausen and Anton Bovier, editors.
\newblock {\em Spin glasses}, volume 1900 of {\em Lecture Notes in
  Mathematics}.
\newblock Springer, Berlin, 2007.

\bibitem{bovier2006statistical}
Anton Bovier.
\newblock {\em Statistical mechanics of disordered systems: a mathematical
  perspective}, volume~18 of {\em Cambridge Series in Statistical and
  Probabilistic Mathematics}.
\newblock Cambridge University Press, Cambridge, 2006.

\bibitem{bovier1998mathematical}
Anton Bovier and Pierre Picco, editors.
\newblock {\em Mathematical aspects of spin glasses and neural networks},
  volume~41 of {\em Progress in Probability}.
\newblock Birkh\"{a}user Boston, Inc., Boston, MA, 1998.

\bibitem{bra83}
J.~G. Brankov and V.~A. Zagrebnov.
\newblock On the description of the phase transition in the
  {H}usimi-{T}emperley model.
\newblock {\em J. Phys. A}, 16(10):2217--2224, 1983.

\bibitem{charbonneau2023spin}
Patrick Charbonneau, Enzo Marinari, Marc Mézard, Giorgio Parisi, Federico
  Ricci-Tersenghi, Gabriele Sicuro, and Francesco Zamponi.
\newblock {\em Spin glass theory and far beyond}.
\newblock World Scientific, 2023.

\bibitem{chen2024on}
Hong-Bin Chen.
\newblock On free energy of non-convex multi-species spin glasses.
\newblock {Preprint, arXiv:2411.13342}.

\bibitem{chen2025convex}
Hong-Bin Chen, Victor Issa, and Jean-Christophe Mourrat.
\newblock The convex structure of the {P}arisi formula for multi-species spin
  glasses.
\newblock {Preprint, arXiv:2508.06397}.

\bibitem{chen2024free}
Hong-Bin Chen and Jean-Christophe Mourrat.
\newblock On the free energy of vector spin glasses with nonconvex
  interactions.
\newblock {\em Probab. Math. Phys.}, 6(1):1--80, 2025.

\bibitem{chen2022hamilton}
Hong-Bin Chen and Jiaming Xia.
\newblock {Hamilton-Jacobi equations from mean-field spin glasses}.
\newblock {\em Probab. Theory Related Fields}, 192(3):803--873, 2025.

\bibitem{chen2022statistical}
Hongbin Chen, Jean-Christophe Mourrat, and Jiaming Xia.
\newblock Statistical inference of finite-rank tensors.
\newblock {\em Ann. H. Lebesgue}, 5:1161--1189, 2022.

\bibitem{chen2018tap}
Wei-Kuo Chen and Dmitry Panchenko.
\newblock On the {TAP} free energy in the mixed {$p$}-spin models.
\newblock {\em Comm. Math. Phys.}, 362(1):219--252, 2018.

\bibitem{chen2023generalized}
Wei-Kuo Chen, Dmitry Panchenko, and Eliran Subag.
\newblock Generalized {TAP} free energy.
\newblock {\em Comm. Pure Appl. Math.}, 76(7):1329--1415, 2023.

\bibitem{coja2018}
Amin Coja-Oghlan.
\newblock Phase transitions in discrete structures.
\newblock In {\em European {C}ongress of {M}athematics}, pages 599--618. Eur.
  Math. Soc., Z\"{u}rich, 2018.

\bibitem{contucci2013perspectives}
Pierluigi Contucci and Cristian Giardin\`a.
\newblock {\em Perspectives on spin glasses}.
\newblock Cambridge University Press, Cambridge, 2013.

\bibitem{dartois2024injective}
Stephane Dartois and Benjamin McKenna.
\newblock Injective norm of real and complex random tensors {I}: From spin
  glasses to geometric entanglement.
\newblock {Preprint, arXiv:2404.03627}.

\bibitem{dominicis2006random}
Cirano De~Dominicis and Irene Giardina.
\newblock {\em Random fields and spin glasses: a field theory approach}.
\newblock Cambridge University Press, New York, 2006.

\bibitem{ding2015proof}
Jian Ding, Allan Sly, and Nike Sun.
\newblock Proof of the satisfiability conjecture for large k.
\newblock In {\em Proceedings of the forty-seventh annual ACM symposium on
  Theory of computing}, pages 59--68, 2015.

\bibitem{ding2016maximum}
Jian Ding, Allan Sly, and Nike Sun.
\newblock Maximum independent sets on random regular graphs.
\newblock {\em Acta Math.}, 217(2):263--340, 2016.

\bibitem{ding2016satisfiability}
Jian Ding, Allan Sly, and Nike Sun.
\newblock Satisfiability threshold for random regular {NAE}-{SAT}.
\newblock {\em Comm. Math. Phys.}, 341(2):435--489, 2016.

\bibitem{sparse_PDE}
Tomas Dominguez and Jean-Christophe Mourrat.
\newblock Infinite-dimensional {H}amilton-{J}acobi equations for statistical
  inference on sparse graphs.
\newblock {\em SIAM J. Math. Anal.}, 56(4):4530--4593, 2024.

\bibitem{sparse_prob}
Tomas Dominguez and Jean-Christophe Mourrat.
\newblock Mutual information for the sparse stochastic block model.
\newblock {\em Ann. Probab.}, 52(2):434--501, 2024.

\bibitem{HJbook}
Tomas Dominguez and Jean-Christophe Mourrat.
\newblock {\em Statistical mechanics of mean-field disordered systems: a
  {H}amilton-{J}acobi approach}.
\newblock Zurich Lectures in Advanced Mathematics. European Mathematical
  Society, Z\"{u}rich, 2024.

\bibitem{elalaoui2021optimization}
Ahmed El~Alaoui, Andrea Montanari, and Mark Sellke.
\newblock Optimization of mean-field spin glasses.
\newblock {\em Ann. Probab.}, 49(6):2922--2960, 2021.

\bibitem{gamarnik2021overlap-survey}
David Gamarnik.
\newblock The overlap gap property: a topological barrier to optimizing over
  random structures.
\newblock {\em Proc. Natl. Acad. Sci. USA}, 118(41):e2108492118, 2021.

\bibitem{gamarnik2025turing}
David Gamarnik.
\newblock Turing in the shadows of {N}obel and {A}bel: an algorithmic story
  behind two recent prizes.
\newblock {\em Notices Amer. Math. Soc.}, 72(5):485--493, 2025.

\bibitem{gamarnik2021overlap-paper}
David Gamarnik and Aukosh Jagannath.
\newblock The overlap gap property and approximate message passing algorithms
  for {$p$}-spin models.
\newblock {\em Ann. Probab.}, 49(1):180--205, 2021.

\bibitem{gamarnik2022disordered}
David Gamarnik, Cristopher Moore, and Lenka Zdeborov\'{a}.
\newblock Disordered systems insights on computational hardness.
\newblock {\em J. Stat. Mech. Theory Exp.}, (11):Paper No. 114015, 41, 2022.

\bibitem{guerra2001sum}
Francesco Guerra.
\newblock Sum rules for the free energy in the mean field spin glass model.
\newblock {\em Fields Institute Communications}, 30(11), 2001.

\bibitem{gue03}
Francesco Guerra.
\newblock Broken replica symmetry bounds in the mean field spin glass model.
\newblock {\em Comm. Math. Phys.}, 233(1):1--12, 2003.

\bibitem{hinton2012practical}
Geoffrey~E. Hinton.
\newblock A practical guide to training restricted {Boltzmann} machines.
\newblock In {\em Neural networks: Tricks of the trade}, pages 599--619.
  Springer, 2012.

\bibitem{hinton2025nobel}
Geoffrey~E. Hinton.
\newblock Nobel lecture: Boltzmann machines.
\newblock {\em Reviews of Modern Physics}, 97(3):030502, 2025.

\bibitem{hopfield1982neural}
John~J. Hopfield.
\newblock Neural networks and physical systems with emergent collective
  computational abilities.
\newblock {\em Proc. Natl. Acad. Sci. USA}, 79(8):2554--2558, 1982.

\bibitem{huang2025tight}
Brice Huang and Mark Sellke.
\newblock Tight {L}ipschitz hardness for optimizing mean field spin glasses.
\newblock {\em Comm. Pure Appl. Math.}, 78(1):60--119, 2025.

\bibitem{huang2023algorithmic}
Brice Huang and Mark Sellke.
\newblock Algorithmic threshold for multi-species spherical spin glasses.
\newblock \noop{2023}{Preprint, arXiv:2303.12172}.

\bibitem{jekel2025potential}
David Jekel, Juspreet~Singh Sandhu, and Jonathan Shi.
\newblock Potential {H}essian ascent: the {S}herrington-{K}irkpatrick model.
\newblock In {\em Proceedings of the 2025 {A}nnual {ACM}-{SIAM} {S}ymposium on
  {D}iscrete {A}lgorithms ({SODA})}, pages 5307--5387. SIAM, Philadelphia, PA,
  2025.

\bibitem{mutual_information}
Anastasia Kireeva and Jean-Christophe Mourrat.
\newblock Breakdown of a concavity property of mutual information for
  non-{G}aussian channels.
\newblock {\em Inf. Inference}, 13(2):Paper No. iaae008, 21, 2024.

\bibitem{krzakala2007gibbs}
Florent Krzaka{\l}a, Andrea Montanari, Federico Ricci-Tersenghi, Guilhem
  Semerjian, and Lenka Zdeborov{\'a}.
\newblock Gibbs states and the set of solutions of random constraint
  satisfaction problems.
\newblock {\em Proc. Natl. Acad. Sci. USA}, 104(25):10318--10323, 2007.

\bibitem{little1974existence}
William~A Little.
\newblock The existence of persistent states in the brain.
\newblock {\em Mathematical biosciences}, 19(1-2):101--120, 1974.

\bibitem{mezmon}
Marc M\'{e}zard and Andrea Montanari.
\newblock {\em Information, physics, and computation}.
\newblock Oxford Graduate Texts. Oxford University Press, Oxford, 2009.

\bibitem{mpvbook}
Marc M\'{e}zard, Giorgio Parisi, and Miguel Virasoro.
\newblock {\em Spin glass theory and beyond}, volume~9 of {\em World Scientific
  Lecture Notes in Physics}.
\newblock World Scientific Publishing Co., Inc., Teaneck, NJ, 1987.

\bibitem{mezard2002analytic}
Marc M{\'e}zard, Giorgio Parisi, and Riccardo Zecchina.
\newblock Analytic and algorithmic solution of random satisfiability problems.
\newblock {\em Science}, 297(5582):812--815, 2002.

\bibitem{monasson1999determining}
R{\'e}mi Monasson, Riccardo Zecchina, Scott Kirkpatrick, Bart Selman, and
  Lidror Troyansky.
\newblock Determining computational complexity from characteristic ‘phase
  transitions’.
\newblock {\em Nature}, 400(6740):133--137, 1999.

\bibitem{montanari2021optimization}
Andrea Montanari.
\newblock Optimization of the {Sherrington-Kirkpatrick Hamiltonian}.
\newblock {\em SIAM Journal on Computing}, (0):FOCS19--1, 2021.

\bibitem{mourrat2020nonconvex}
Jean-Christophe Mourrat.
\newblock Nonconvex interactions in mean-field spin glasses.
\newblock {\em Probab. Math. Phys.}, 2(2):281--339, 2021.

\bibitem{mourrat2022parisi}
Jean-Christophe Mourrat.
\newblock The {P}arisi formula is a {H}amilton-{J}acobi equation in
  {W}asserstein space.
\newblock {\em Canad. J. Math.}, 74(3):607--629, 2022.

\bibitem{mourrat2023free}
Jean-Christophe Mourrat.
\newblock Free energy upper bound for mean-field vector spin glasses.
\newblock {\em Ann. Inst. Henri Poincar\'{e} Probab. Stat.}, 59(3):1143--1182,
  2023.

\bibitem{mourrat2024informal}
Jean-Christophe Mourrat.
\newblock An informal introduction to the {P}arisi formula.
\newblock {Preprint, arXiv:2410.12364}.

\bibitem{mourrat2024uninverting}
Jean-Christophe Mourrat.
\newblock Un-inverting the {P}arisi formula.
\newblock {\em Ann. Inst. H. Poincar\'{e} Probab. Statist.}, to appear.

\bibitem{mourrat2020extending}
Jean-Christophe Mourrat and Dmitry Panchenko.
\newblock Extending the {P}arisi formula along a {H}amilton-{J}acobi equation.
\newblock {\em Electron. J. Probab.}, 25:Paper No. 23, 17, 2020.

\bibitem{mulet2002coloring}
Roberto Mulet, Andrea Pagnani, Martin Weigt, and Riccardo Zecchina.
\newblock Coloring random graphs.
\newblock {\em Phys. Rev. Lett.}, 89(26):268701, 2002.

\bibitem{new86}
Charles Newman.
\newblock Percolation theory: A selective survey of rigorous results.
\newblock In {\em Advances in multiphase flow and related problems}. SIAM,
  1986.

\bibitem{nishimori2001statistical}
Hidetoshi Nishimori.
\newblock {\em Statistical physics of spin glasses and information processing},
  volume 111 of {\em International Series of Monographs on Physics}.
\newblock Oxford University Press, New York, 2001.

\bibitem{opper2001advanced}
Manfred Opper and David Saad, editors.
\newblock {\em Advanced mean field methods}.
\newblock Neural Information Processing Series. MIT Press, Cambridge, MA, 2001.

\bibitem{pan.aom}
Dmitry Panchenko.
\newblock The {P}arisi ultrametricity conjecture.
\newblock {\em Ann. of Math. (2)}, 177(1):383--393, 2013.

\bibitem{pan}
Dmitry Panchenko.
\newblock {\em The {S}herrington-{K}irkpatrick model}.
\newblock Springer Monographs in Mathematics. Springer, New York, 2013.

\bibitem{pan.multi}
Dmitry Panchenko.
\newblock The free energy in a multi-species {S}herrington-{K}irkpatrick model.
\newblock {\em Ann. Probab.}, 43(6):3494--3513, 2015.

\bibitem{pan.potts}
Dmitry Panchenko.
\newblock Free energy in the {P}otts spin glass.
\newblock {\em Ann. Probab.}, 46(2):829--864, 2018.

\bibitem{pan.vec}
Dmitry Panchenko.
\newblock Free energy in the mixed {$p$}-spin models with vector spins.
\newblock {\em Ann. Probab.}, 46(2):865--896, \noop{2019}2018.

\bibitem{parisi1979infinite}
Giorgio Parisi.
\newblock Infinite number of order parameters for spin-glasses.
\newblock {\em Phys. Rev. Lett.}, 43(23):1754, 1979.

\bibitem{parisi1980order}
Giorgio Parisi.
\newblock The order parameter for spin glasses: a function on the interval 0-1.
\newblock {\em J. Phys. A}, 13(3):1101, 1980.

\bibitem{parisi1980sequence}
Giorgio Parisi.
\newblock A sequence of approximated solutions to the {S-K} model for spin
  glasses.
\newblock {\em J. Phys. A}, 13(4):L115--L121, 1980.

\bibitem{parisi1983order}
Giorgio Parisi.
\newblock Order parameter for spin-glasses.
\newblock {\em Phys. Rev. Lett.}, 50(24):1946, 1983.

\bibitem{parisi2023nobel}
Giorgio Parisi.
\newblock Nobel lecture: Multiple equilibria.
\newblock {\em Reviews of Modern Physics}, 95(3):030501, 2023.

\bibitem{richardson2008modern}
Tom Richardson and Ruediger Urbanke.
\newblock {\em Modern coding theory}.
\newblock Cambridge university press, 2008.

\bibitem{schoenberg1942positive}
I.~J. Schoenberg.
\newblock Positive definite functions on spheres.
\newblock {\em Duke Math. J.}, 9:96--108, 1942.

\bibitem{sellke2021optimizing}
Mark Sellke.
\newblock Optimizing mean field spin glasses with external field.
\newblock {\em Electron. J. Probab.}, 29:Paper No. 4, 47, 2024.

\bibitem{sherrington1975solvable}
David Sherrington and Scott Kirkpatrick.
\newblock Solvable model of a spin-glass.
\newblock {\em Phys. Rev. Lett.}, 35(26):1792, 1975.

\bibitem{smolensky1986information}
Paul Smolensky.
\newblock Information processing in dynamical systems: foundations of harmony
  theory.
\newblock In {\em Parallel distributed processing: explorations in the
  microstructure of cognition}, volume~1, pages 194--281. MIT Press, 1986.

\bibitem{stein2013spin}
Daniel~L. Stein and Charles~M. Newman.
\newblock {\em Spin glasses and complexity}.
\newblock Primers in Complex Systems. Princeton University Press, Princeton,
  NJ, 2013.

\bibitem{subag2021following}
Eliran Subag.
\newblock Following the ground states of full-{RSB} spherical spin glasses.
\newblock {\em Comm. Pure Appl. Math.}, 74(5):1021--1044, 2021.

\bibitem{subag2023tap2}
Eliran Subag.
\newblock T{AP} approach for multispecies spherical spin glasses {II}: {T}he
  free energy of the pure models.
\newblock {\em Ann. Probab.}, 51(3):1004--1024, 2023.

\bibitem{subag2025tap1}
Eliran Subag.
\newblock T{AP} approach for multi-species spherical spin glasses {I}:
  {G}eneral theory.
\newblock {\em Electron. J. Probab.}, 30:Paper No. 87, 32, 2025.

\bibitem{Tpaper}
Michel Talagrand.
\newblock The {P}arisi formula.
\newblock {\em Ann. of Math. (2)}, 163(1):221--263, 2006.

\bibitem{Tbook1}
Michel Talagrand.
\newblock {\em Mean field models for spin glasses. {V}olume {I}}, volume~54 of
  {\em Ergebnisse der Mathematik und ihrer Grenzgebiete}.
\newblock Springer-Verlag, Berlin, 2011.

\bibitem{Tbook2}
Michel Talagrand.
\newblock {\em Mean field models for spin glasses. {V}olume {II}}, volume~55 of
  {\em Ergebnisse der Mathematik und ihrer Grenzgebiete}.
\newblock Springer, Heidelberg, 2011.

\bibitem{thouless1977solution}
David~J Thouless, Philip~W Anderson, and Robert~G Palmer.
\newblock Solution of'solvable model of a spin glass'.
\newblock {\em Philosophical Magazine}, 35(3):593--601, 1977.

\bibitem{tubiana2018restricted}
J{\'e}r{\^o}me Tubiana.
\newblock {\em {Restricted Boltzmann machines: from compositional
  representations to protein sequence analysis}}.
\newblock PhD thesis, {Universit{\'e} Paris Sciences et Lettres}, 2018.

\bibitem{tubiana2017emergence}
J{\'e}r{\^o}me Tubiana and R{\'e}mi Monasson.
\newblock Emergence of compositional representations in restricted boltzmann
  machines.
\newblock {\em Phys. Rev. Lett.}, 118(13):138301, 2017.

\bibitem{zdekrz}
Lenka Zdeborov{\'a} and Florent Krzakala.
\newblock Statistical physics of inference: Thresholds and algorithms.
\newblock {\em Advances in Physics}, 65(5):453--552, 2016.

\end{thebibliography}

\end{document}